\documentclass{amsart}

\usepackage{amsmath, amsthm}
\usepackage{amssymb}
\usepackage{amsfonts}
\usepackage{latexsym}
\usepackage{mathrsfs}
\usepackage{bm}
\usepackage{graphicx}
\usepackage{leftidx}
\usepackage{mathtools}

\theoremstyle{plain}
\newtheorem{theorem}{Theorem}[section]
\newtheorem*{theorem*}{Theorem}
\newtheorem{lem}[theorem]{Lemma}

\newtheorem{cor}{Corollary}

\theoremstyle{remark}
\newtheorem{Rem}{Remark}[section]

\newtheorem{exm}{Example}

\newcommand{\csf}{curve shortening flow}

\newcommand{\Rm}{Riemannian manifold}

\numberwithin{equation}{section}

\newcommand{\be}{\begin{equation}}
\newcommand{\ene}{\end{equation}}
\newcommand{\br}{\begin{Rem}}
\newcommand{\er}{\end{Rem}}
\newcommand{\bl}{\begin{lem}}
\newcommand{\el}{\end{lem}}
\newcommand{\bd}{\begin{Def}}
\newcommand{\ed}{\end{Def}}
\newcommand{\ben}{\begin{enumerate}}
\newcommand{\een}{\end{enumerate}}
\newcommand{\bp}{\begin{proof}}
\newcommand{\ep}{\end{proof}}
\newcommand{\bpo}{\begin{pro}}
\newcommand{\epo}{\end{pro}}
\newcommand{\beq}{\begin{equation*}}
\newcommand{\eeq}{\end{equation*}}
\newcommand{\bear}{\begin{eqnarray*}}
\newcommand{\eear}{\end{eqnarray*}}
\newcommand{\bt}{\begin{theorem}}
\newcommand{\et}{\end{theorem}}
\newcommand{\bst}{\begin{split}}
\newcommand{\est}{\end{split}}
\newcommand{\bal}{\begin{aligned}}
\newcommand{\eal}{\end{aligned}}

\newcommand{\F}[2]{\frac{#1}{#2}}
\renewcommand{\P}{\partial}
\renewcommand{\c}{\mathcal{C}}

\renewcommand{\H}{\mathcal{H}}
\newcommand{\N}{\mathcal{N}}

\newcommand{\PT}[1]{\frac{\partial}{\partial #1}}

\newcommand{\PQ}[2]{\frac{\partial #1}{\partial #2}}

\newcommand{\ct}{\c_{t}}
\newcommand{\co}{\c_{0}}
\newcommand{\Lo}{\F{d}{dt}-\triangle^{\ct}}
\newcommand{\vs}{ds_{t}}
\newcommand{\e}{\epsilon}

\newcommand{\Pt}{\F{\P}{\P t}}
\def\XXint#1#2#3{{\setbox0=\hbox{$#1{#2#3}{\int}$}
    \vcenter{\hbox{$#2#3$}}\kern-.5\wd0}}

\makeatletter
\def\@citestyle{\m@th\upshape\mdseries}
\def\citeform#1{{\bfseries#1}}
\def\@cite#1#2{{%
  \@citestyle[\citeform{#1}\if@tempswa, #2\fi]}}
\@ifundefined{cite }{%
  \expandafter\let\csname cite \endcsname\cite
  \edef\cite{\@nx\protect\@xp\@nx\csname cite \endcsname}%
}{}
\makeatother
\begin{document}
\title{The Curve shortening Flow with Parallel 1-form}

\author{Hengyu Zhou}
\address{The Graduate Center,\\ The City University of New York,\\ 365 fifth Ave.,\\ New York, NY 10016, USA}
\email{hzhou@gc.cuny.edu}

\begin{abstract}
     Let $M$ be a closed Riemannian manifold with a parallel 1-form $\Omega$. We prove two theorems about the curve shortening flow in $M$. One is that the {\csf} $\ct$ in $M$ exists for all $t$ in $[0, \infty)$, if it satisfies $\Omega(T)\geq 0$ on the initial curve $\co$. Here $T$ is the unit tangent vector on $\co$. The other one is about the convergence. It says that in a closed {\Rm} $\tilde{M}$, assume the curve shortening flow $\ct$ exists for all $t\in[0,\infty)$  and its length converges to a positive limit, then  $ \lim\limits _{t\rightarrow\infty}max_{\ct}|\nabla^{m}A|^{2}=0$ for all $m=0,1,\cdots$. Here $A$ denotes the second fundamental form of $\ct$ in $\tilde{M}$.
\end{abstract}
\date{\today}

\maketitle
\section*{Introduction}

  We use $\nabla$ denote the Levi-Civita connection of a Riemannian manifold . A 1-form $\Omega$ is called parallel if $\nabla\Omega=0$. In this paper we consider the long time existence and the convergence of the {\csf} in the Riemannian manifold with arbitrary codimension. More precisely, the {\csf} $\ct$  is given as the solution of the follow equation:
     \be\label{eq:gmcf}
    \left\{\bal
     &\F{\P\c_{t}}{\P t}=\vec{H}\\
             &\c(x ,0)=\c_{0}(x);
             \eal\right.
        \ene
    Here $\vec{H}=\nabla_{T}T$ is the mean curvature vector of $\ct$ in the ambient Riemannian manifold.

    The evolution of closed curves under \eqref{eq:gmcf} has received considerable study. For example, Gage (\cite{Gag84}), Gage-Hamilton (\cite{GH86}) and
   Grayson (\cite{Gra87}), considered the evolution of the convex curve flowed by the mean curvature in the plane. Furthermore, Grayson (\cite{Gra89}) generalized the results into the case on the embedded closed curve of Riemannian surfaces. Huisken (\cite{Hui94}) and Andrews-Bryan (\cite{BB11}) used the technique of distance comparison to investigate the {\csf} on the plane. For  the mean curvature flow of the submanifolds with codimension $k\geq2$, Andrews-Baker (\cite{BC10}) proposed a machinery about the evolution of the second fundamental form $A$ with arbitrary codimension. They also obtained the convergence of the mean curvature flow for the submanifold  pinched to sphere in Euclidean space. More recently, with an integral curvature condition (\cite{LX11}) and a pinched condition in hyperbolic space form (\cite{LX12}) Liu, Xu, Ye and Zhao also proved two extension theorems for the mean curvature flow of the submanifold with high codimension. Smoczyk's survey (\cite{Sm12}) is a good reference.

   In a different direction, Wang (\cite{Wan02}) initialized a way to use parallel n-forms $\Omega$, $n\geq 2$, to investigate the graphic mean curvature flow for arbitrary codimension ($\geq 2$) in product manifolds. With the evolution equation of $*\Omega$, ($*$ is the Hodge dual), Wang obtained its long-time existence and an estimate of its second fundamental form along the mean curvature flow. Moreover, Wang and his coauthors generalized the results of graphic mean curvature flow into other settings (\cite{MW11},\cite{TW04},\cite{Wan01},\cite{SW02},).
   Both of (\cite{BC10}) and (\cite{Wan02}) did not  consider the case of the curve shortening flow, which is the topic of this paper.

     There are two special properties about the {\csf} whose codimension in its ambient manifold is greater than 1.  First the embeddedness of the evolving curve is not necessarily preserved (\cite{Alt91}). Therefore, we only consider the case of the immersed curve. On the other hand, $|\vec{H}|^{2}=|A|^{2}$. Based on this fact, we can overcome the barrier that the curve shortening flow has arbitrary codimension. This makes many computations possible, especially in Section 2.

    Now let us state two main theorems:

    \bt \label{TH1}
     Let $M$ be a closed Riemannian manifold, and $\Omega$ be a parallel 1-form in $M$. If $\co$ is a closed curve satisfying $\Omega(T)>0$ on $\co$. Here $T$ is the unit tangent vector on $\co$. Then the curve shortening flow $\ct$ in \eqref{eq:gmcf} exists for all $t$ in $[0,\infty)$.
     \et

     \begin{cor} \label{Mcora}
      Let $(N,g)$ be a closed manifold, $du^{2}$ is the canonical metric on  $S^{1}$. Assume $M=S^{1}\times N$ is the product manifold with the product metric $\bar{g}=g+du^{2}$. Let $\Omega=du$, which is a parallel 1-form in $M$. $\co$ is a closed curve in $M$. If the unit tangent vector $T$ of $\co$ satisfies $\Omega(T)>0$, the {\csf} $\ct$ in \eqref{eq:gmcf} exists for all $t$ in $[0, \infty)$.
     \end{cor}

     \br
  Corollary \ref{Mcora} can be viewed the 1-dimension version of Theorem A in (\cite{Wan02}). The existence result here with $\Omega(T)>0$ is stronger than that of graphical mean curvature flow with $*\Omega>\frac{1}{2}$ in (\cite{Wan02}). It also can be compared with the results of Tsui-Wang (\cite{TW04}) if we choose $N=S^{n}$.
     \er
Now we consider the convergence of the {\csf} in any closed {\Rm}. Notice that we do not require there is a parallel 1-form in the ambient manifold. The convergence of the {\csf} in the $C^{\infty}$ sense can be stated as follows.
   \bt \label{TH2}
     Let $\tilde{M}$ be a closed manifold. If the {\csf} $\ct$ in $\tilde{M}$  exists for all $t\in[0, \infty)$ and the length of $\ct$  converges to a positive number, $\lim\limits_{t\rightarrow \infty}max_{\c_{t}}|\nabla^{m}A|=0$ for all $m=0,1,\cdots$.
   \et
   \br
     The definition of $\nabla^{m}A$ is very subtle. In (\cite{Hui84},\cite{Hui86}), Huisken's definitions about $\nabla^{m} A$ only work in the hypersurface's case. Andrews-Bakers (\cite{BC10}) gave a rigorous definition of $\nabla^{m} A$ for any mean curvature flow with arbitrary codimension in Riemannian manifold. The basic idea is  that $|\nabla^{m}A|$  can contain all information about the mean curvature flow. Namely, when $|\nabla^{m}A|$ goes to 0 when $t$ goes to infinity, the mean curvature flow $F_{t}$ will converge to a geodesic submanifold in some natural sense. For example, one of these cases is that $F_{t}$ has long time existence and always stay in a compact region of its ambient Riemannian manifold. We will state the definition of $\nabla^{m}A$ in Section 2.
   \er
   \br
    The condition that $M$ and $\tilde{M}$ are closed is not essential. We take them only for the sake of the exposition. In Theorem \ref{TH1}, what we really need is that for all finite time $t$, the curve shortening flow $\ct$ is in a compact region in $M$, which can depend on time $t$. In Theorem \ref{TH2}, we can only require that all covariant derivatives of the Riemann tensor in $\tilde{M}$ are uniformly bounded.
   \er

 \subsection*{Outline of this paper}

 In section 1, we give two examples to illustrate our motivations of using  parallel form to investigate the curve shortening flow. In Section 2, we give the definition of $\nabla^{m}A$. Then we derive the evolution equations related to the second fundamental form $A$ (lemma \ref{EVA} and lemma \ref{EVnorm}). With those equations, we derive the evolution equation of $\Omega(T)$ for any 1-form $\Omega$ along the curve shortening flow (lemma \ref{EVform}). In particular, when $\Omega$ is a parallel 1-form, this evolution equation states that the lower positive bound of $\Omega(T)$ can be preserved along the curve shortening flow. In section 3, we prove Theorem \ref{TH1} by using the results from section 2. In section 4, with Grayson's idea (\cite{Gra89}), we prove Theorem \ref{TH2}.

\subsection*{Notation} We collect some geometric quantities and facts used later.
 \begin{enumerate}
 \item   $\c(u,t)$ denotes the solution of \eqref{eq:gmcf}.  $\F{\P}{\P t}$ and $\F{\P}{\P u}$ are the abbreviations of $\c_{*}(\F{\P}{\P t})$ and $\c_{*}(\F{\P}{\P u})$ respectively.
\item $\nabla$ is the Levi-Civita connection, $du^{2}$ is the canonical metric on $S^{1}$.
\item  Let $R$ denote Riemann curvature tensor, given by  $R(X,Y)=\nabla_{Y}\nabla_{X}-\nabla_{X}\nabla_{Y}+\nabla_{[X,Y]}$.
\item Let $s$ be the arc-length parameter of a curve $\c$. Its unit tangent vector $T$ is equal to $ \F{\P}{\P s}=\F{\c_{*}(\F{\P}{\P u})}{<\c_{*}(\F{\P}{\P u}),\c_{*}(\F{\P}{\P u})>^{\F{1}{2}}}.$
 \item  $\vec{H}=\F{\P}{\P t}=A(T,T)=\nabla_{T}T$. $|A|^{2}=|\vec{H}|^{2}$.
 \item  For a curve $\c$, $|\vec{H}|^{2}=|A|^{2}=<\nabla_{T}T,\nabla_{T}T>$. All curves in this paper are closed.
 \end{enumerate}

\subsection*{Acknowledgements}  The author would like to thank Professor Zeno Huang for suggesting the problem, and for many
helpful discussions and encouragement, and Professor Yunping Jiang for his interest and suggestions. This work is completed while the author is
partially supported, as a student associate, by a CIRG award $\# 1861$. The author is also indebted to Xiangwen Zhang, Longzhi Lin and Andy Huang. These discussions with them greatly help to simplify the computations in this paper.

   \section{The Examples about Parallel form}
   In this section we use two examples to illustrate the motivations that  we choose parallel 1-form $\Omega$ to investigate the {\csf}. The first one is the graphic mean curvature flow in the product manifolds (\cite{Wan02}). Based on this, we can use Wang's idea (\cite{Wan02}) about parallel form to propose a problem about  the open curve shortening flow in $R^{1+p}$.

    \begin{exm}[graphic mean curvature flow] In (\cite{Wan02}), Mu-Tao Wang  initialized a way to use parallel form to deform a map along the mean curvature flow in product manifolds. The basic setting is as followed:
Let $N_{1}$, $N_{2}$ be closed manifolds with constant sectional curvatures, dimension $\geq 2$.  $f$ is a smooth map from $N_{1}$ to $N_{2}$. Let $F(p,t): N_{1}\times[0,q)\rightarrow N_{1}\times N_{2}$ satisfies the following equation.
      \be\label{GQ}
     \left\{\bal
      \frac{\P F_{t}}{\P t}&=\vec{H}\\
              F(p ,0)&=(p,f(p));
     \eal
      \right.
      \ene
      Let $\Omega$ be the volume form of $N_{1}$, which is parallel in the product manifold $N_{1}\times N_{2}$. Together with some technical conditions about the sectional curvature of $N_{1}$ and $N_{2}$, Wang (\cite{Wan02}) obtained that if $*\Omega>\frac{1}{2}$ at the initial manifold, the evolution equations of  $*\Omega$ indicates  the long time existence of the solution in \eqref{GQ}.
    \end{exm}

  \begin{exm}[open curve shortening flow]
  In this example, we propose a question about open curve shortening flow in $R^{1+p}$.
 $\co$ in $R^{1+p}$ is defined by
   \be \label{GC}
   \c_{0}:R^{1}\rightarrow R^{1+p}\quad \mbox{and}\qquad \c_{0}(x)=(x,f_{1}(x),\cdots,f_{p}(x));
   \ene
   Let $|Df|^{2}=\sum_{i=1,\cdots,p}(f_{i}^{'}(x))^{2}$. Then,
   $$
   T=\left(\frac{\partial}{\partial x^{1}}+\sum_{i=1,
   \cdots,p}f_{i}^{'}(x)\frac{\partial}{\partial x^{i+1}}\right)\frac{1}{\sqrt{1+|Df|^{2}}};$$

    Let $\Omega=dx_{1}$. Again $\Omega$ is a parallel 1-form. On $\c_{0}$,  $\Omega(T)=\frac{1}{\sqrt {1+|Df|^{2}}}$. From lemma \ref{EP}, if $\ct$ exists for $t\in [0,q)$, then $\Omega(T)$ on $\ct$ satisfies
        $$ \F{\P }{\P t}\Omega(T)=\triangle^{\ct}\Omega (T)+|A|^{2}\Omega(T)$$
    Here $A$ is its second fundamental form.
     Recall that with $*\Omega>\F{1}{2}$ in the initial data, Wang (\cite{Wan02}) answered the long time existence of graphic mean curvature flow positively. Naturally if we assume that $\co$ satisfies  $\Omega(T)\geq\delta>0$, does the {\csf}
        $\ct$ in \eqref{eq:gmcf} exist for all $t\in [0,\infty)$? If it exists, can we give any description about the convergence of this curve shortening flow in $R^{1+p}$?

  \end{exm}

   \br
    The above examples reflect our basic perspectives about the parallel form and the mean curvature flow.  \eqref{EP} and \eqref{GF} are some possible evidences to support these connections between them. \eqref{EP} can be viewed as a generalization of Wang's graphical mean curvature flow (\cite{Wan02}) in the curve's case. \eqref{GF} indicates that $\mu=\F{1}{\Omega(T)}$ be thought as the gradient function (\cite{EH89}) in the case of the curve shortening flow. There are no any oblivious connections between (\cite{Wan02}) and (\cite{EH89}). For the curve shortening flow, however, \eqref{GF} and \eqref{EP} are two different forms of the same evolution equation.
     \er

  \section{Evolution Equations}

   We derive the evolution equations of the curve shortening flow. Generally, the forms of those equations are very complicated for the mean curvature flow whose dimension and codimension are both $> 1$.  (\cite{BC10}, \cite{Wan02})

  \subsection{Definition of $\nabla^{m}A$}

   We briefly state the definition of $\nabla^{m}A$ which can be found in (\cite{BC10}). Assume $M$ is a Riemannian manifold with the Levi-Civita connection $\bar{\nabla}$, $N$ is an immersed submanifold of $M$. Suppose $I$ is a real interval, then the tangent space $T(N\times I)$ splits into a direct product $\H\otimes R\P t$, where $\H=\{u\in T(N\times I); dt(u)=0\}$ is the "spatial" tangent bundle.

   We consider a smooth map $F:N\times I\rightarrow M$ which is a time-dependent immersion, i.e.,for each fixed $t\in I$, $F(\,,t):N\rightarrow M$ is an immersion. $F^{*}TM$ is a vector bundle over $N\times I$. We can define a metric $g_{F}$ and the connection $^{F}\nabla$ on $F^{*}TM$ by the pull-back from $(\bar{\nabla},\bar{g})$ on $M$. Let $\N$ be the orthogonal complement of $\H$ in $F^{*}TM$. We denote the $\pi,\pi^{\bot}$ be the orthogonal projections from $F^{*}TM$ onto $\H$ and $\N$ respectively. The connections $^{\H}\nabla$ and $^{N}\nabla$ are given by.
   $$
   \bal
        ^{\H}\nabla&=\pi\circ^{F}\nabla\circ F_{*};\\
         ^{\N}\nabla&=\pi^{\bot}\circ^{F}\nabla\circ F_{*};
   \eal
   $$
   For a tensor $K\in \Gamma(\H^{*}\otimes \H^{*}\otimes N)$. We have to define $\nabla^{m}K\in \Gamma(\otimes^{m+2}\H^{*}\otimes N)$ by induction on $m$. Let $u_{0}, u_{i},i=1,2,3,\cdots,m+1$, then $\nabla^{m}K$ and $\nabla_{\F{\P}{\P t}}\nabla^{m}K$ are given by
   $$
   \bal
     \nabla_{u_{0}}&(\nabla^{m-1}K)(u_{1},\cdots,u_{m+1})=^{\N}\nabla_{u_{0}}(\nabla^{m-1}K(u_{1},\cdots,u_{m+1}))-\sum\limits_{i=1}^{m+1}\nabla^{m-1}K(\cdots ,^{\H}\nabla_{u_{0}}u_{i},\cdots);\\
     \nabla_{\F{\P}{\P t}}&(\nabla^{m}K)(u_{0},u_{1},\cdots,u_{m+1})=^{\N}\nabla_{\F{\P}{\P t}}(\nabla^{m}K(u_{0},\cdots,u_{m+1}))-\sum\limits_{k=0}^{m+1}\nabla^{m}K(\cdots,^{\H}\nabla_{\F{\P}{\P t}}u_{k},\cdots);
     \eal
     $$

\subsection{The Evolution Equations for the Curve Shortening Flow}
 For the curve shortening flow $\ct$ in \eqref{eq:gmcf}, we define a new tensor $\tilde{R}\in \Gamma(\H^{*}\otimes\H^{*}\otimes\N)$ by $\tilde{R}(T,T)=R(T,\vec{H})T$. We also have the definition of $\nabla^{m}\tilde{R}$. Since $A\in \Gamma(\H^{*}\otimes\H^{*}\otimes\N)$, with the definitions in the previous subsection, we can obtain a theorem about the evolution equation of $A$ from equation (18) in (\cite{BC10}).
 \bt [Andrews-Baker]
 For the curve shortening flow $\ct$, the evolution equation of $A(T,T)$ is given by
     \be\label{eSFF}
   \nabla_{\Pt}A(T,T)=\nabla_{T}\nabla_{T} A(T,T)+|A|^{2}\vec{H}+\tilde{R}(T,T);
     \ene
 \et

Recall that  $S*T$ (\cite{Hui84}) means any linear combination of tensors formed by contraction on $S$ and $T$. \eqref{eSFF} can be rewritten as $\nabla_{\Pt}A=\nabla^{2} A+ \sum A*A*A+\tilde{R}(T,T)$. This form can be generalized into any $m$-th covariant derivative $A$ for the curve shortening flow $\ct$.
\bl  \label{EVA}
The evolution equation of $m$-th covariant derivative $A$ is of the form.
\be
\nabla_{\Pt}\nabla^{m}A=\nabla^{m+2}A+\sum\limits_{i+j+k=m}\nabla^{i}A*\nabla^{j}A*\nabla^{k}A+\nabla^{m}\tilde{R}(T,T);
\ene
\el
\bp
  We prove the lemma through the induction on $m$. The case $m=0$ is given by \eqref{eSFF}. Now suppose the results hold up to $k\leq m-1$.
  In (\cite{Hui84}), the time derivative of the Christoffel symbols $\Gamma_{jk}^{i}$ has the form $A*\nabla A$. Therefore we obtain,

  \begin{align*}
  \nabla_{\Pt}\nabla^{m}A&=\nabla(\nabla_{\Pt}\nabla^{m-1}A)+\nabla^{m-1}A*A*\nabla A;\\
                         &=\nabla\{\nabla^{m+1}A)+\sum\limits_{i+j+k=m-1}\nabla^{i}A*\nabla^{j}A*\nabla^{k}A+\nabla^{m-1}\tilde{R}(T,T)\}\\
                         &+\nabla^{m-1}A*A*\nabla A;\\
                         &=\nabla^{m+2}A+\sum\limits_{i+j+k=m}\nabla^{i}A*\nabla^{j}A*\nabla^{k}A+\nabla^{m}\tilde{R}(T,T);
   \end{align*}
                         \ep
Moreover, we obtain the evolution equations of $|\nabla^{m}A|^{2},m=0,1,\cdots,$.
\bl \label{EVnorm}
The evolution equations of $|\nabla^{m}A|^{2}$ are given by
\begin{gather}
 \F{\P}{\P t}|A|^{2}=\triangle^{\ct}|A|^{2}-2|\nabla A|^{2}+2|A|^{4}+2R(T,\vec{H},T,\vec{H});\label{EQSF}\\
\begin{split}
 \F{\P}{\P t}|\nabla^{m}A|^{2}&=\triangle^{\ct}|\nabla^{m}A|^{2}-2|\nabla ^{m+1}A|^{2}+\sum\limits_{i+j+k=m}\nabla^{i}A*\nabla^{j}A*\nabla^{k}A*\nabla^{m}A\\
                      &+2<\nabla^{m}\tilde{R}(T,T),\nabla^{m}A>\label{EQSFF}\quad\text{for}\quad m\geq 1;
\end{split}
\end{gather}
\el
\bp
   We differentiate $|\nabla^{m}A|^{2}$ with respect to $t$. For $m=0$,
   $$
   \bal
   \F{\P}{\P t}|A|^{2}&=2<\nabla_{\Pt}A(T,T),A(T,T)>;\\
                      &=2<\nabla^{2}A(T,T)+|A|^{2}\vec{H}+R(T,\vec{H})T,A(T,T)>;\\
    \text{Since $A(T,T)=\vec{H}$}\\
                      &=\triangle^{\ct}|A|^{2}-2|\nabla A|^{2}+2|A|^{4}+2R(T,\vec{H},T,\vec{H});
   \eal
   $$
   For $m\geq 1$.
   $$
   \bal
   \F{\P}{\P t}|\nabla^{m}A|^{2}&=2<\nabla_{\Pt}\nabla^{m}A,\nabla^{m}A>;\\
                      &=2<\nabla^{m+2}A+\sum\limits_{i+j+k=m}\nabla^{i}A*\nabla^{j}A*\nabla^{k}A+\nabla^{m}\tilde{R}(T,T),\nabla^{m}A>;\\
                      &=\triangle^{\ct}|\nabla^{m}A|^{2}-2|\nabla ^{m+1}A|^{2}+\sum\limits_{i+j+k=m}\nabla^{i}A*\nabla^{j}A*\nabla^{k}A*\nabla^{m}A+\\
                      &+2<\nabla^{m}\tilde{R}(T,T),\nabla^{m}A>;
   \eal
   $$

\ep

          Now we derive the evolution equation of $\Omega(T)$ along the {\csf} $\ct$.
         \bl \label{EVform}
         Let $\Omega$ be a 1-form in a {\Rm} $\tilde{M}$, $s$ be the arc-length parameter of $\ct$. Let $\triangle^{\ct}=\F{\P^{2}}{\P s^{2}}$ be the Laplace operator on $\ct$, and $T=\F{\P}{\P s}$ be the unit tangent vector. Thus we have
         \be \label{EP1}
             \F{\P }{\P t}\Omega(T)=\triangle^{\ct}\Omega (T)+|A|^{2}\Omega(T)-\nabla^{2}\Omega(T,T,T)
             -2\nabla\Omega(T,\nabla_{T}T);
         \ene
         In particular, if $\Omega$ is parallel in $\tilde{M}$, we obtain
          \be \label{EP}
             \F{\P }{\P t}\Omega(T)=\triangle^{\ct}\Omega (T)+|A|^{2}\Omega(T);
         \ene
          \el
         \bp
  Recall $T=\F{C_{*}(\F{\P}{\P u})}{<C_{*}(\F{\P}{\P u}),C_{*}(\F{\P}{\P u})>^{\F{1}{2}}}$ and $\F{\P}{\P t}=\vec{H}$.  And
     \be
     \begin{split}
         [\F{\P}{\P t},T]&=\F{[C_{*}(\F{\P}{\P t}),C_{*}(\F{\P}{\P u})]}{<C_{*}(\F{\P}{\P u}),C_{*}(\F{\P}{\P u})>^{\F{1}{2}}}+C_{*}(\F{\P}{\P u})\F{\P}{\P t}(\F{1}{<C_{*}(\F{\P}{\P u}),C_{*}(\F{\P}{\P u})>^{\F{1}{2}}});\\
              &=-C_{*}(\F{\P}{\P u})\F{<\nabla_{C_{*}(\F{\P}{\P t})}(C_{*}(\F{\P}{\P u})),C_{*}(\F{\P}{\P u})>}{<C_{*}(\F{\P}{\P u}),C_{*}(\F{\P}{\P u})>^{\F{3}{2}}};\\
              &=-C_{*}(\F{\P}{\P u})\F{<\nabla_{C_{*}(\F{\P}{\P u})}(\vec{H}),C_{*}(\F{\P}{\P u})>}{<C_{*}(\F{\P}{\P u}),C_{*}(\F{\P}{\P u})>^{\F{3}{2}}};\\
              &=|A|^{2}T;
     \end{split}
     \ene
Therefore, the relation of $\nabla_{\F{\P}{\P t}}T$ and $\nabla_{T}\F{\P}{\P t}$ is given by
 \be
   \nabla_{\F{\P}{\P t}}T-\nabla_{T}\F{\P}{\P t}=|A|^{2}T;
 \ene

     \be  \label{lap}
            \begin{split}
         \triangle^{\ct}\Omega(T)&=T\big(\Omega(\nabla_{T}T)+\nabla\Omega(T,T)\big);\\
                                    &=\Omega(\nabla_{T}\F{\P}{\P t})+\nabla\Omega(\F{\P}{\P t},T)
                                    +\nabla^{2}\Omega(T,T,T)+2\nabla\Omega(T,\nabla_{T}T);
         \end{split}
         \ene
         The $t$-derivative of $\Omega(T)$

            \be \label{eq:Ev2}
            \begin{split}
              \F{d}{dt}\Omega(T)&=\Omega(\nabla_{\frac{\partial }{\partial t}}T)+\nabla\Omega(\F{\P }{\P t},T); \\
                          &=\Omega(\nabla_{T}\F{\P }{\P t})+|A|^{2}\Omega(T)+\nabla\Omega(\F{\P }{\P t},T);\quad
              \end{split}
            \ene
           Therefore, \eqref{lap} and  \eqref{eq:Ev2} lead to the lemma.
        \ep
     \br
       It's easily checked $min_{\ct}\Omega$ is a Lipschitz function of $t$, and $min_{\ct}\Omega(T)$ is differentiable with respect to $t$ almost everywhere. When $\Omega$ ia a parallel 1-form in $\tilde{M}$,
       $\Omega(T)>0$ for $t\in [0,q_{0})$. For $ \triangle^{\ct}\Omega(T)\geq 0$ for the points which attain the minimal value of $\Omega(T)$, for a.e $t\in [0,q_{0})$ we have the following:
            $$\F{d}{d t}min_{\ct}\Omega(T)\geq |A|^{2}min_{\ct}\Omega(T)\geq 0;$$
        This implies that $min_{\ct}\Omega(T)$ is nondecreasing along the flow. Then $\Omega(T)\geq \delta_{0}$ will be preserved whenever the curve shortening flow exists. 
     \er

  \section{Long time existence}

  With the evolution results in the previous section, we prove Theorem \ref{TH1}. The short time existence of $\ct$ is from the short time existence of solution of the quasilinear parabolic equation for closed initial data. It's well known that if the mean curvature flow $F_{t}$ of the hypersurface in Euclidean space exists for only finite time interval $[0,q)$, $max_{F_{t}}|A|\rightarrow\infty$ when $t$ approach to $q$ (\cite{Hui84}). For the submanifold in Euclidean space with  high codimension, such kind of lemma is proved by Andrews-Baker (\cite{BC10}). For the curve shortening flow's case, we give the proof of the following lemma only for the sake of completeness.
\bl\label{UnboundedSFF}
Let $\tilde{M}$ be a {\Rm}.
 If the curve shortening flow $\ct$ in $\tilde{M}$ exists for $t$ in the maximal finite interval $[0,q)$, then $ max_{\ct}|A|^{2}\rightarrow \infty \ \text{as}\ t \rightarrow q$.

\el
\bp
If the lemma is false, there exists a constant $C_{4}<\infty$  such that $ max_{\ct}|A|\leq C_{4}$ for $t\in[0,q)$. It follows that for all $u\in S^{1}$ and $t_{1},t_{2}\in[0,q)$.
     $$
     \bal
      dist(\c_{t_{1}}(u),\c_{t_{2}}(u))&\leq |\int\limits_{t_{1}}^{t_{2}}\vec{H}dt|\leq C_{4}|t_{1}-t_{2}|;\\
        \F{d}{dt}\int_{\ct}ds_{t}&\geq -C_{4}\int_{\ct}ds_{t};
        \eal
        $$
 Therefore, $\ct$ converges uniformly to some continuous limit $\c_{q}(u)$ and the length $l_{\ct}\geq\delta>0$ for $t\in [0,q)$.
  We want to show that $\c_{q}(u)$ actually represents a smooth limit curve. Then we can extend the flow $\ct$ over time $q$. It's a contradiction to the maximal finite interval in the assumption.

 By \eqref{EQSFF}, we are led to
 \be \label{e2}
 \begin{split}
          \F{\P}{\P t}|\nabla^{m}A|^{2}&\leq \triangle^{\ct}|\nabla^{m}A|^{2}-2|\nabla^{m+1}A|^{2}
                +C_{m}\sum_{i+j+k=m}|\nabla^{i}A||\nabla^{j}A||\nabla^{k}A||\nabla^{m}A|\\
                &+D_{m}\sum_{j\leq m}|\nabla^{j}A||\nabla^{m}A|+\tilde{C}_{m}|\nabla^{m}A|;
        \end{split}
 \ene

   By Cauchy inequality, we have the follows.
     \be\label{e3}
     \bal |\nabla^{m}A|^{2}&\leq min_{\ct}|\nabla^{m}A|^{2}\vs+\int_{\ct}\F{\P}{\P s}|\nabla^{m}A|^{2}\vs;\\
                          &\leq \F{1}{l_{\ct}}\int |\nabla^{m}A|^{2}\vs+2 \int_{\ct}<\nabla^{m}A,\nabla^{m+1}A> \vs;\\
                           &\leq (\F{1}{\delta}+2)\int |\nabla^{m}A|^{2}\vs+2 \int_{\ct}|\nabla^{m+1}A|^{2} \vs;
             \eal
        \ene

   If $\ct$ satisfies $|\nabla^{m}A|\leq C_{m}$ for all $t<q$, then $\c_{q}$ is a smooth curve. From our definition about $\nabla^{m}A$, the proof of this fact is classical but a little tedious. Therefore we omit it here. For the Euclidean case, please refer to the proof of Theorem 3 in (\cite{BC10}).

   By \eqref{e3}, we have to prove that $\int_{\ct}|\nabla^{m}A|^{2}ds_{t}\leq C_{m}$ for all $m$.
   We argue $|\int_{\ct}\nabla^{m}A\vs|\leq C_{m}$ by induction on $m$. If it holds up for $k\leq m-1$. Then by \eqref{e3} for $k\leq m-2$, $max_{\ct}|\nabla^{k}A|\leq C_{k}$ for $t\in [0,q)$. \eqref{e2} and \eqref{e3} imply that
     $$ \F{\P}{\P t}\int_{ct}|\nabla^{m}A|^{2}ds_{t}\leq \tilde{C}(m)\big(\int_{ct}|\nabla^{m}A|^{2}ds_{t}+(\int_{ct}|\nabla^{m}A|^{2}ds_{t})^{\F{1}{2}}\big);
$$
  With the above inequality, we obtain
       $$
       (\int_{\ct}|\nabla^{m}A|^{2}ds_{t})^{\F{1}{2}}+1\leq e^{\F{\tilde{C}(m)}{2}t}((\int_{\co}|\nabla^{m}A|^{2}ds_{0})^{\F{1}{2}}+1);
       $$
  Then we conclude that for all $m=1,\cdots, max_{\ct}|\nabla^{m}A|\leq C_{m}$. $C_{q}(u)$ is the smooth limit curve of $\ct$ for $t<q$. This will lead to the contradiction, since we suppose $[0,q)$ is the maximal time interval for the existence of $\ct$.
\ep

Now let $\mu=\F{1}{\Omega(T)}$. The following lemma tells us $\mu$ can be viewed as the gradient function $\nu$ in (\cite{EH89}). They have the same form of evolution equation.
    \bl \label{L1}
    Let $M$ be a closed manifold, $\Omega$ be a parallel 1-form in $M$, $s$ be the arc-length parameter of $\ct$.
      If the {\csf} $\ct$ exists for $t\in [0,q)$ and $\co$ satisfies $\Omega(T)>0$,
     then $\mu$ satisfies
     \be \label{GF}
     (\Lo)\mu=-|A|^{2}\mu-2\mu^{-1}|\F{\P\mu}{\P s}|^{2};
     \ene
 \el
 \bp  $$\bal
 (\Lo)\mu&=-\F{1}{\Omega(T)^{2}}(\Lo)\Omega(T)-2\Omega(T)|\PQ{\Omega(T)}{s}|;\\
     &=-\F{1}{\Omega(T)^{2}}(|A|^{2}\Omega(T))-2\Omega(T)|\PT{s}\F{1}{\Omega(T)}|^{2};\\
     &=-|A|^{2}\mu-2\mu^{-1}|\F{\P\mu}{\P s}|^{2};
 \eal
 $$
 \ep
   \bl
    Use the assumption in lemma \ref{L1}. Let $C_{0}=max_{x\in M}R$. If the curve shortening flow $\ct$ exists for $t\in [0, q)$, then
   \be
   (\Lo)|A|^{2}\mu^{2}\leq -2\mu^{-1}\PQ{\mu}{s}\PQ{|A|^{2}\mu^{2}}{s}+2C_{0}|A|^{2}\mu^{2};
   \ene
   \el
   \bp
   $$
  (\Lo)|A|^{2}=-2|\nabla A|^{2}+2|A|^{4}+2R(T,\vec{H},T,\vec{H});
 $$
 and together with the identity
     $$(\Lo)\mu^{2}=-2|A|^{2}\mu^{2}-6|\PQ{\mu}{s}|^{2}$$
   yields
        \be \label{3.5}
        (\Lo)|A|^{2}\mu^{2}\leq -2|\nabla A |^{2}\mu^{2}-6|\PQ{\mu}{s}|^{2}|A|^{2}-2\PQ{|A|^{2}}{s}\PQ{\mu^{2}}{s}+2R(T,\vec{H},T,\vec{H})\mu^{2};
        \ene
       and
      \be\label{3.6}
      \bal
      -2\PQ{|A|^{2}}{s}.\PQ{\mu^{2}}{s}&\leq-\PQ{|A|^{2}}{s}\PQ{\mu^{2}}{s}-4\mu |A|\PQ{|A|}{s}\PQ{\mu}{s};\\
      & \leq-\mu^{-2}\PQ{\mu^{2}}{s}\PQ{|A|^{2}\mu^{2}}{s}+4|\PQ{
      \mu}{s}|^{2}|A|^{2}-4\mu |A|\PQ{|A|}{s}\PQ{\mu}{s}\\
      &\leq  -2\mu^{-1}\PQ{\mu}{s}\PQ{|A|^{2}\mu^{2}}{s}+2|\PQ{|A|}{s}|^{2}\mu^{2}+6|\PQ{\mu}{s}|^{2}|A|^{2};
     \eal
      \ene
      Here we use the fact $\F{\P |A|}{\P s}\leq |\nabla A|$.
   For $M$ is closed, $|\vec{H}|^{2}=|A|^{2}$, $|R(T,\vec{H},T,\vec{H})|\leq C_{0}|A|^{2}$.
     By \eqref{3.5} and \eqref{3.6}, we obtain
     \be\label{MPC}
     (\Lo)|A|^{2}\mu^{2}\leq -2\mu^{-1}\PQ{\mu}{s}\PQ{|A|^{2}\mu^{2}}{s}+2C_{0}|A|^{2}\mu^{2};
     \ene
   \ep
   \begin{cor}\label{C1}
   Under the assumption of Theorem \ref{TH1}. Then
   $max_{\ct}|A|\leq C(q,C_{0})$ when $t<q$. Here $C(q,C_{0})$ is a constant only dependent on $q$ and $C_{0}$.
   \end{cor}
   \bp  Because $M$ is closed, for any unit vector field $S$ in $TM$, $|\Omega(S)|\leq C_{M}$.
    $$
   \mu^{-1}\PQ{\mu}{s}=\Omega(T)\PT{s}\F{1}{\Omega(T)}=-\mu\Omega(\vec{H});$$
   Since $|\mu\Omega(\vec{H})|$ is continuous for $(u,t)\in S^{1}\times [0,q)$,
      the maximal principle of parabolic equation in \eqref{MPC} gives ,
     \be\label{thenormofsecond}
      \F{\P}{\P t}max_{\ct}|A|^{2}\mu^{2}\leq 2C_{0}max_{\ct}|A|^{2}\mu^{2};
      \ene
    From $\mu\geq \F{1}{C_{M}}$ and \eqref{thenormofsecond},  we conclude the corollary.
   \ep
\subsubsection*{The proof of Theorem 1}
Now we prepare everything which we need to prove Theorem \ref{TH1}. From corollary \ref{C1}, for all finite time interval $[0,q)$, $max_{\ct}|A|$ is always uniformly bounded (depending on $q$). If $\ct$ exists only for a finite interval, Lemma \ref{UnboundedSFF} will lead to a contraction.  Then Theorem 1 is followed.

         \section{The Convergence}

     This section is devoted to prove Theorem \ref{TH2}. Recall  Theorem \ref{TH2}  concludes that the limit of $max_{\c_{t}}|\nabla^{m}A|$ is $0$ for all $m=0,1,\cdots$  if  the {\csf} $\ct$ in the closed manifold $\tilde{M}$ exists for all $t\in[0, \infty)$ and the length of $\ct$  converges to a positive number. The idea of its proof originated from Section 7 in (\cite{Gra89}).  We suppose the length of $\ct$ satisfies $l_{\ct}\geq l_{\infty}>0$. First, let's prove a technique lemma.
        \bl \label{T3}
        For all $m\geq 0$,

        \be \label{Ta}
        |\nabla^{m}A|^{2}\leq a \int|\nabla^{m}A|^{2}ds_{t}+2\int |\nabla^{m+1}A|^{2}ds_{t};
        \ene
        Here $a=\F{1}{l_{\infty}}+2$.  For $m\geq 1$,
         \be \label{Tb}(\int |\nabla^{m}A|^{2}\vs)^{2}\leq \int |\nabla^{m-1}A|^{2}\vs \int |\nabla^{m+1}A|^{2}\vs;
          \ene

      \el
      \bp
         For $|\nabla^{m}A|^{2}$,
         $$ \bal |\nabla^{m}A|^{2}&\leq min_{\ct}|\nabla^{m}A|^{2}\vs+\int_{\ct}\F{\P}{\P s}|\nabla^{m}A|^{2}\vs;\\
                          &\leq \F{1}{l_{\ct}}\int |\nabla^{m}A|^{2}\vs+2 \int_{\ct}<\nabla^{m}A,\nabla^{m+1}A> \vs;\\
                           &\leq (\F{1}{l_{\ct}}+2)\int |\nabla^{m}A|^{2}\vs+2 \int_{\ct}|\nabla^{m+1}A|^{2} \vs;
             \eal
        $$
        Integrating $\int|\nabla^{m}A|^{2}\vs$ by parts,
          $$\int|\nabla^{m}A|^{2}\vs=\int \F{\P}{\P s}<\nabla^{m-1}A,\nabla^{m}A>\vs-\int <\nabla^{m-1}A,\nabla^{m+1}A>\vs;$$
         The first term is 0. Using the Cauchy inequality we obtain \eqref{Tb}.
      \ep

    \subsubsection*{The proof of Theorem \ref{TH2}}  Our strategy is to prove for all $m=0,1,\cdots$, $\int|\nabla^{m}A|^{2}\vs$ converges to 0 as $t\rightarrow \infty$. First, we prove $m=0$ case, then argue by induction.

      From the defintion of the curve shortening flow, we get the following two facts.
          \begin{gather}
          \F{d}{dt}\int_{\ct} ds_{t}=-\int|\vec{H}|^{2}ds_{t}=-\int_{\ct}|A|^{2}ds_{t}; \\
        l_{\ct}-l_{\co}=\int\limits_{0}^{t}\F{d}{dt}\int_{\ct} ds_{t}dt=-\int\limits_{0}^{t}\int_{\ct}|A|^{2}ds_{t}dt<\infty;
        \end{gather}

   There exists a measure zero set $E\,\text{in}\,[0,\infty)$. Let $E_{n}=[n,\infty)\cap E^{c}$. We get
   \be\label{AMD}
   max_{t\in E_{n}}\int_{\ct}|A|^{2}ds_{t}\rightarrow 0\quad \text{as}\quad n\rightarrow\infty;
   \ene
From \eqref{EQSF}, we differentiate $\int_{\ct}|A|^{2}ds_{t}$ with respect to $t$.
        \be
        \bal
         \F{d}{d t}\int_{\ct}|A|^{2}ds_{t}&\leq -2\int|\nabla A|^{2}ds_{t}+2\int_{\ct}|A|^{4}ds_{t}+D_{0}\int|A|^{2}ds_{t};\\
         &\leq -2\int|\nabla A|^{2}ds_{t}+ (a\int_{\ct}|A|^{2}ds_{t}+2\int_{\ct}|\nabla A|^{2}ds_{t})\int|A|^{2}ds_{t}\\
         &+D_{0}\int|A|^{2}ds_{t};\\
         &\leq-(2-2\int_{\ct}|A|^{2}ds_{t})\int_{\ct}|\nabla A|^{2}ds_{t}+\int|A|^{2}ds_{t}(D_{0}+a\int|A|^{2}ds_{t});
        \eal
        \ene
 Here we use \eqref{Ta}. Let $n$ be a sufficiently large number. We can assume $(2-2\int_{\ct}|A|^{2}ds_{t})\geq 1$ and $D_{0}+a\int|A|^{2}ds_{t}\leq D_{0}+1$ for all $t\in E_{n}$. Furthermore,

      \be\label{epgrowth}
       \F{d}{d t}\int_{\ct}|A|^{2}ds_{t}\leq (D_{0}+1)\int_{\ct}|A|^{2}ds_{t}
       \ene
 For both sides of \eqref{epgrowth} are continuous with respect to $t$, \eqref{epgrowth} holds up for $t\in [n,\infty)$. $\int_{\ct}|A|^{2}ds_{t}$ increases at most with exponential growth when $n$ is sufficiently large.
 By \eqref{AMD}, we obtain $\lim_{t\rightarrow \infty}\int |A|^{2}ds_{t}=0$.

We notice that $\tilde{M}$ is closed, all $|\nabla^{m}R|\leq C_{m}\,\text{for all}\, m$. Together with \eqref{EQSFF}, we have the following inequality for $m=1,\cdots$.
        \begin{multline}\label{CDSFF}
          \F{\P}{\P t}|\nabla^{m}A|^{2}\leq \triangle^{\ct}|\nabla^{m}A|^{2}-2|\nabla^{m+1}A|^{2}\\
                +\tilde{C}_{m}\sum_{i+j+k=m}|\nabla^{i}A||\nabla^{j}A||\nabla^{k}A||\nabla^{m}A|+D_{m}\sum_{j\leq m}|\nabla^{j}A||\nabla^{m}A|+\tilde{C}_{m}|\nabla^{m}A|;
        \end{multline}
Here all constants in this proof are only depending on $C_{m}$.

When $m\geq 1$, we argue by the induction on $m$. Suppose $\lim\limits_{t\rightarrow \infty}\int_{\ct}|\nabla^{k}A|^{2}ds_{t}= 0$ for $k\leq m-1$. Then by \eqref{Ta}, for $k\leq m-2$, $\lim\limits _{t\rightarrow \infty}max_{\ct}|\nabla^{k}A|^{2}=0$.

For any $\e>0$, let $t_{0}$ be a positive number sufficiently large, such that for $t\geq t_{0}$ and $k\leq m-2$, $\int_{\ct}|\nabla^{m-1}A|^{2}ds_{t}\leq \e$ and $max_{\ct}|\nabla^{k}A|^{2}\leq \e$. By \eqref{CDSFF}, we obtain the following estimates.
$$
\bal
\F{\P}{\P t}|\nabla^{m}A|^{2}&\leq \triangle^{\ct}|\nabla^{m}A|^{2}-2|\nabla^{m+1}A|^{2}+\tilde{C}_{m}(\e)(|\nabla^{m}A|^{2}\\
&+|\nabla^{m}A||\nabla^{m-1}A|)+ \tilde{D}_{m}|\nabla^{m}A|;\\
\F{d}{dt}\int_{\ct}|\nabla^{m}A|^{2}ds_{t}&\leq -2\int_{\ct}|\nabla^{m+1}A|^{2}ds_{t}+\tilde{C}_{m}(\e)(\int _{\ct}|\nabla^{m}A|^{2}ds_{t}\\
&+\int_{\ct}|\nabla^{m}A||\nabla^{m-1}A|ds_{t})+\tilde{D}_{m}\int_{\ct}|\nabla^{m}A|ds_{t};
\eal
$$
 Here $\tilde{C}_{m}(\e)$ goes to 0 as $\e$ converges to 0.

  Let $C>3$ be a constant determined later. If $\int_{\ct}|\nabla^{m}A|^{2}ds_{t}\geq C\int_{\ct}|\nabla^{m-1}A|^{2}ds_{t}$, the following estimates are given by \eqref{Tb}.
 \be
 \bal
 \int_{\ct}|\nabla^{m}A|^{2}ds_{t}&\leq\F{1}{C}\int_{\ct}|\nabla^{m+1}A|^{2}ds_{t};\\
 \int_{\ct}|\nabla^{m-1}A|^{2}ds_{t}&\leq\F{1}{C^{2}}\int_{\ct}|\nabla^{m+1}A|^{2}ds_{t};
 \eal
 \ene
 As a result,
  \be
    \bal
    \int_{\ct}|\nabla^{m}A||\nabla^{m-1}A|ds_{t}&\leq\F{1}{\sqrt{C^{3}}}\int_{\ct}|\nabla^{m+1}A|^{2}ds_{t};\\
    \int_{\ct}|\nabla^{m}A|ds_{t}&\leq\F{\sqrt{l_{\co}}}{\sqrt{C}}(\int_{\ct}|\nabla^{m+1}A|^{2}ds_{t})^{\F{1}{2}};
    \eal
  \ene
 The $t$-derivative of $\int_{\ct}|\nabla^{m}A|^{2}ds_{t}$ can be written as followed:
 \be\label{Mfinish}
 \bal
 \F{d}{dt}\int_{\ct}|\nabla^{m}A|^{2}ds_{t}&\leq -(2-\tilde{C}_{m}(\e)(\F{1}{\sqrt{C^{3}}}+\F{1}{C}))\int_{\ct}|\nabla^{m+1}A|^{2}ds_{t}\\
 &+\F{\tilde{D}_{m}}{\sqrt{C}}(\int_{\ct}|\nabla^{m+1}A|^{2}ds_{t})^{\F{1}{2}};
 \eal
\ene
Let $C$ be sufficiently large such that $(2-\tilde{C}_{m}(\e)(\F{1}{\sqrt{C^{3}}}+\F{1}{C}))\geq 1$ and $\F{\tilde{D}_{m}}{\sqrt{C}}\leq 2\e$. \eqref{Mfinish} becomes
   \be\label{finish}
 \F{d}{dt}\int_{\ct}|\nabla^{m}A|^{2}ds_{t}\leq -\int_{\ct}|\nabla^{m+1}A|^{2}ds_{t})^{\F{1}{2}} (\int_{\ct}|\nabla^{m+1}A|^{2}ds_{t})^{\F{1}{2}}- 2\e);
\ene
 Now we can conclude Theorem \ref{TH2} as follows.
\begin{enumerate}
\item
If $\int_{\ct}|\nabla^{m}A|^{2}ds_{t}\geq C\int_{\ct}|\nabla^{m-1}A|^{2}ds_{t}$,
\begin{enumerate}

\item If $(\int_{\ct}|\nabla^{m}A|^{2}ds_{t})^{\F{1}{2}}\geq 3\e$, we have $(\int_{\ct}|\nabla^{m+1}A|^{2}ds_{t})^{\F{1}{2}}\geq 3\e$.  From \eqref{finish}, then $\F{d}{dt}\int_{\ct}|\nabla^{m}A|^{2}ds_{t}\leq -2\e^{2}$. $\int_{\ct}|\nabla^{m}A|^{2}ds_{t}$ will decrease until
    $$(\int_{\ct}|\nabla^{m}A|^{2}ds_{t})\leq (3\e)^{2}$$
 If we can find a $t_{0}$ such that
   $$(\int_{\ct}|\nabla^{m}A|^{2}ds_{t})^{\F{1}{2}}\leq 3\e;$$
    then for all $t>t_{0}$, $\int_{\ct}|\nabla^{m}A|^{2}ds_{t}\leq (3\e)^{2}$.
\item  If $(\int_{\ct}|\nabla^{m}A|^{2}ds_{t})^{\F{1}{2}}\leq 3\e$, the conclusion of theorem \ref{TH2} is obviously true.
\end{enumerate}

\item If $\int_{\ct}|\nabla^{m}A|^{2}ds_{t}< C\int_{\ct}|\nabla^{m-1}A|^{2}ds_{t}$, we have $\int_{\ct}|\nabla^{m}A|^{2}ds_{t}\leq C\e$.

 \end{enumerate}

 In a word, finally $\int_{\ct}|\nabla^{m}A|^{2}ds_{t}$ converges to 0 when $t\rightarrow \infty$. By \eqref{Ta}, $max_{\ct}|\nabla^{m}A|$ converges to 0 for all $m$ as $t$ goes to $\infty$.

\bibliographystyle{alpha}	
%\bibliography{isombedding}
\bibliography{ref}
  \end{document}